\documentclass[11pt]{amsart}
\usepackage{amssymb, latexsym}
\usepackage{graphicx}
\theoremstyle{plain}
\newtheorem{theorem}{Theorem}

\newtheorem*{thm-cheb}{Theorem (Chebyshev)}

\newtheorem{proposition}{Proposition}

\newtheorem*{2'}{Theorem 2'}
\newtheorem*{3'}{Theorem 3'}

\theoremstyle{remark}

\newtheorem*{Remark 1}{Remark 1}
\newtheorem*{Remark 2}{Remark 2}
\newtheorem*{Remark 3}{Remark 3}
\newtheorem*{Remark 4}{Remark 4}

\numberwithin{equation}{section}

\begin{document}


\title[Permutations avoiding a pattern under   Mallows distributions]
{Permutations avoiding a  pattern   of length three under   Mallows distributions }
\author{Ross G. Pinsky}


\dedicatory{Dedicated to the memory of Dima Ioffe (1963-2020)}

\address{Department of Mathematics\\
Technion---Israel Institute of Technology\\
Haifa, 32000\\ Israel}
\email{ pinsky@math.technion.ac.il}

\urladdr{http://www.math.technion.ac.il/~pinsky/}

\date{}

\begin{abstract}
We consider permutations avoiding a  pattern of length three under the family of Mallows distributions. In particular, for any pattern $\tau\in S_3\setminus\{321\}$,
we obtain rather precise results on the asymptotic probability as $n\to\infty$ that
a permutation $\sigma\in S_n$ under the Mallows distribution with parameter $q\in(0,1)$ avoids
the pattern.
By a  duality between the parameters $q$ and $\frac1q$,
we  also obtain rather precise results on  the above probability for $q>1$ and any pattern $\tau\in S_3\setminus \{123\}$.

\end{abstract}
\maketitle

\section{Introduction and Statement of Results}\label{intro}
We recall the definition of pattern avoidance for permutations. Let $S_n$ denote the set of permutations of $[n]:=\{1,\cdots, n\}$.
For $\sigma\in S_n$, we write $\sigma=\sigma_1\sigma_2\cdots\sigma_n$, where $\sigma(i)=\sigma_i$.  If $\tau\in S_m$, where $2\le m\le n$,
then we say that $\sigma$ contains $\tau$ as a pattern if there exists a subsequence $1\le i_1<i_2<\cdots<i_m\le n$ such
that for all $1\le j,k\le m$, the inequality $\sigma_{i_j}<\sigma_{i_k}$ holds if and only if the inequality $\tau_j<\tau_k$ holds.
If $\sigma$ does not contain $\tau$, then we say that $\sigma$ \it avoids\rm\ $\tau$.
Thus, for example, if $\tau=123\in S_3$, then  $\sigma=53412\in S_5$ avoids $\tau$ but $\sigma=51324\in S_5$ does not avoid $\tau$; indeed the pattern $\tau$ appears twice: $\sigma_2\sigma_3\sigma_5=134$ and $\sigma_2\sigma_4\sigma_5=124$.
Denote by $S_n(\tau)$  the set of permutations in $S_n$ that avoid $\tau$.
At a couple of  points in this paper it will be useful to have the extended definition
$S_n(\tau)=S_n$, if $\tau\in S_m$ with $m>n$.

It is well-known  that
$|S_n(\tau)|=C_n$,  for all six permutations $\tau\in S_3$, where $C_n=\frac{\binom{2n}n}{n+1}$ is the $n$th
Catalan number \cite{B}. Consequently, under the uniformly random probability measure $P_n$ on $S_n$, one has
\begin{equation}\label{asympunif}
P_n(S_n(\tau))\sim
\frac{(4e)^n}{\sqrt2\thinspace \pi n^{n+2}},\ \text{for all}\ \tau\in S_3.
\end{equation}

More generally, the celebrated Stanley-Wilf conjecture, now proven \cite{MT}, states that for any $\tau\in S_m$, $m\ge2$, there exists a constant $C=C(\tau)$ such that
$C^n$ constitutes an upper bound on the growth rate of $|S_n(\tau)|$ as $n\to\infty$.
Consequently, for any $\tau\in S_m$, $m\ge2$, there exists a constant $C_1=C_1(\tau)$ such that
$P_n(S_n(\tau))\le (\frac{C_1}n)^n$.

In this paper, we investigate the probability of avoiding a  pattern of length three under the family of Mallows distributions.
For each $q>0$, the Mallows distribution with parameter $q$ is the probability measure $P_n^q$ on $S_n$ defined by
\begin{equation}\label{mallowsdef}
P_n^q(\sigma)=\frac{q^{\text{inv}(\sigma)}}{Z_n(q)},
\end{equation}
where $\text{inv}(\sigma)$ is the number of inversions in $\sigma$, and $Z_n(q)$ is the normalization constant \cite[Corollary 1.3.13]{S}, given by
\begin{equation}\label{Zn}
Z_n(q)=\prod_{k=1}^n\frac{1-q^k}{1-q}.
\end{equation}
Thus, for $q\in(0,1)$, the distribution favors permutations with few inversions, while for $q>1$, the distribution favors permutations with many inversions.
Recall that the reverse  of a permutation $\sigma=\sigma_1\cdots\sigma_n$ is the permutation $\sigma^{\text{rev}}:=\sigma_n\cdots\sigma_1$.
The Mallows distributions satisfy the following duality between $q>1$ and $q<1$:
$$
P_n^q(\sigma)=P_n^{\frac1q}(\sigma^{\text{rev}}),\ \text{for}\ q>0, \sigma\in S_n\ \text{and}\ n=1,2,\cdots.
$$
Consequently, it suffices to restrict our study to $q\in (0,1)$, which we do from here on.
For the above results and more on the Mallows distribution, see for example, \cite{R}, \cite{GP} \cite{BP}.
We are unaware of other papers treating the probability of avoiding a pattern of length three under the Mallows
distributions. However, the paper \cite{CD} treats the probability of avoiding a \it consecutive\rm\ pattern under the Mallows distribution.

For several of the proofs in this paper, we will use the following so-called online construction of a random permutation in $S_n$ distributed according to the Mallows distribution with parameter $q$. By ``online'' we mean that
the random permutation is constructed in $n$ steps, with one number being added to the permutation at each step.
Let $\{X_j\}_{j=1}^n$ be independent random variables with $X_j$ distributed as a geometric random variable with parameter $1-q$ and truncated at $j-1$; that is,
\begin{equation}\label{Xdist}
P(X_j=m)=\frac{(1-q)q^m}{1-q^j},\ m=0,\cdots, j-1.
\end{equation}
Consider a horizontal line on which to place the numbers in $[n]$.
We begin by placing down the number 1. Then inductively, if  we have already placed down the numbers $1,2,\cdots, j-1$, the number $j$ gets placed down in the position for which there are $X_j$ numbers to its right.
Thus, for example, for $n=4$, if $X_2=1$, $X_3=2$ and $X_4=0$, then we obtain the permutation 3214.
To see that this construction does indeed induce the Mallows distribution with parameter $q$, note
that the number of inversions in the constructed permutation $\sigma$ is $\sum_{j=2}^nX_j$, and thus using \eqref{Xdist},
 $P(X_j=x_j,\ j=2,\cdots, n)=\frac1{Z_n(q)}q^{\sum_{j=2}^nx_j}=\frac{q^{\text{inv}(\sigma)}}{Z_n(q)}$.
\medskip

We  obtain rather precise results for every pattern except 321. The pattern 321 seems to require different techniques and is currently being studied by a post-doctoral student of mine \cite{HP}.
(By duality, for $q>1$, we obtain rather precise results  for every pattern except 123.)
We begin with a very rough result to set the stage.
From the online construction and \eqref{Xdist}, it follows that
\begin{equation}\label{id}
P_n^q(\sigma=id)=P(X_j=0,j=1,\cdots, n)=\frac{(1-q)^n}{\prod_{j=1}^n(1-q^j)},
\end{equation}
where id denotes the identity permutation.
Since the identity permutation avoids all patterns of length three except for the pattern 123,
the above calculation yields the following result.
\begin{proposition}\label{nonidentprop}
For $q\in(0,1)$,
\begin{equation}\label{nonident}
P_n^q(S_n(\tau))> (1-q)^n, \ \text{for all}\ \tau\in S_3\setminus\{123\}\ \text{and all}\ n\ge3.
\end{equation}
\end{proposition}
On the other hand, we will prove the following result regarding the pattern 123.
\begin{theorem}\label{123thm}
For $q\in(0,1)$,
\begin{equation}\label{123}
\lim_{n\to\infty}\big(P_n^q(S_n(123))\big)^\frac1{n^2}= q^{\frac14}.
\end{equation}
\end{theorem}
Proposition \ref{nonidentprop} and Theorem \ref{123thm} along with \eqref{asympunif} show that under the Mallows distribution with $q\in(0,1)$,
$\tau$-avoiding permutations, for $\tau\in S_3\setminus\{123\}$, are overwhelming more likely  than under the uniform distribution, whereas 123-avoiding permutations are overwhelmingly less likely
than under the uniform distribution.

Since the permutation 123 has no inversions, the permutations 132 and 213 have one inversion, the permutations 312 and 231 have two inversions, and the permutation 321 has three inversions,
it is rather natural to  suspect that $P_n^q(S_n(123))<P_n^q(S_n(132))=P_n^q(S_n(213))<P_n^q(S_n(312))=P_n^q(S_n(231))<P_n^q(S_n(321))$, for all $n\ge3$.
We will prove these relations with regard to four  of the six patterns.
\begin{proposition}\label{equals}
Let $q\in(0,1)$. Then
\begin{equation}\label{pairequal}
P_n^q(S_n(312))=P_n^q(S_n(231))>P_n^q(S_n(213))=P_n^q(S_n(132)),\ \text{for all}\ n\ge3.
\end{equation}
\end{proposition}
\bf\noindent Remark.\rm\
If $\sigma$ is distributed as $P_n^q$, then $\sigma^{-1}$ is also distributed as $P_n^q$; indeed this follows from the definition of $P_n^q$ and the fact that the number of inversions in any permutation is equal to the number of inversions in its inverse.
Thus,  since 312 and 231 are inverses of each other, this gives the first equality in \eqref{pairequal}.

Proposition \ref{equals}  actually follows as a corollary of Proposition \ref{basicprop} below. Its proof will be given
after the statement of that proposition.

We have the following proposition.
\begin{proposition}\label{dn1n}
Let $q\in(0,1)$. Then
\begin{equation}\label{dn1nexists}
\lim_{n\to\infty}\big(P_n^q(S_n(\tau))\big)^\frac1n\  \text{exists  for all}\ \tau\in S_3.
\end{equation}
\end{proposition}
\bf \noindent Remark.\rm\ For $\tau=123$, we already know from Theorem \ref{123thm} that the limit in \eqref{dn1nexists} exists and equals 0. From Proposition \ref{nonidentprop}, it follows that the limit  in \eqref{dn1nexists} is positive for
$\tau\in S_3\setminus\{123\}$.
\begin{proof}
For $\sigma\in S_{n_1+n_2}$, $I_1=[n_1]$ and $I_2=[n_1+n_2]\setminus[n_1]$, let $\sigma_{I_i},\ i=1,2$, denote the permutation in $S_{n_i}$ obtained by
the induced relative ordering of $\sigma$ restricted to the domain $I_i$.
(Thus, if $\sigma=32451$,  $n_1=2$ and $n_2=3$, then $\sigma_{I_1}=21$ and $\sigma_{I_2}=231$.)
 It is known \cite[Lemma 2.5 and Corollary 2.7]{BP} that if
$\sigma$ is distributed as $P_{n_1+n_2}^q$, then $\sigma_{I_i}$ is distributed as $P_{n_i}^q$, and $\sigma_{I_1}$ and
$\sigma_{I_2}$ are independent. From this we conclude that
$$
P_{n_1+n_2}^q(S_{n_1+n_2}(\tau))\le P_{n_1}^q(S_{n_1}(\tau))P_{n_2}^q(S_{n_2}(\tau)),\
\text{for any}\ \tau\in S_3.
$$
Thus, $\{\rho_n\}_{n=1}^\infty$ is a sub-additive sequence, where $\rho_n=\log P_n^q(S_n(\tau))$. By the well-known result on sub-additive sequences,
$\lim_{n\to\infty}\frac{\rho_n}n$ exists, and consequently, so does $\lim_{n\to\infty}\big(P_n^q(S_n(\tau))\big)^\frac1n$.
\end{proof}
From Propositions \ref{nonidentprop} and \ref{dn1n}, it follows that
$\lim_{n\to\infty}\big(P_n^q(S_n(\tau))\big)^\frac1n\ge1-q$, for all $\tau\in S_3\setminus\{123\}$. In fact, there is equality
for two choices of $\tau$. We will prove the following theorem.
\begin{theorem}\label{132or213}
$$
\lim_{n\to\infty}\big(P_n^q(S_n(132))\big)^\frac1n=\lim_{n\to\infty}\big(P_n^q(S_n(213))\big)^\frac1n=1-q.
$$
\end{theorem}
The next proposition concerns $P_n^q(S_n(312))$ and $P_n^q(S_n(213))$,
or equivalently by Proposition \ref{equals}, $P_n^q(S_n(231))$ and $P_n^q(S_n(132))$.
Define
\begin{equation}\label{w}
w_n=w_n(q)=\prod_{l=1}^n(1-q^l),n=0,1,\cdots,
\end{equation}
where we use the convention $w_0=1$.
\begin{proposition}\label{basicprop} Let $q\in(0,1)$.

\noindent i.
Define $d_n=d_n(q)=P_n^q(S_n(312))$ or $d_n=d_n(q)=P_n^q(S_n(231))$. Then
\begin{equation}\label{basic1}
d_n=(1-q)\sum_{k=1}^nq^{k-1}\thinspace\frac{w_{k-1}w_{n-k}}{w_n}d_{k-1}d_{n-k},\ n=1,2,\cdots.
\end{equation}
\noindent ii.
 Define $d_n=d_n(q)=P_n^q(S_n(213))$ or $d_n=d_n(q)=P_n^q(S_n(132))$. Then
\begin{equation}\label{basic2}
d_n=(1-q)\sum_{k=1}^nq^{(n-k+1)(k-1)}\thinspace\frac{w_{k-1}w_{n-k}}{w_n}d_{k-1}d_{n-k},\ n=1,2,\cdots.
\end{equation}
\end{proposition}
Proposition \ref{equals} is  a direct corollary of Proposition \ref{basicprop}.

\noindent \it Proof of Proposition \ref{equals}.\rm\
The proposition follows immediately by induction from \eqref{basic1} and \eqref{basic2} along with the fact that
$P_n^q(S_n(\tau))=1$ for $n=1,2$ and  all $\tau\in S_3$.\hfill $\square$
\medskip

We can use \eqref{basic1} to study  $\lim_{n\to\infty}\big(P_n^q(S_n(312))\big)^\frac1n$
(or equivalently, $\lim_{n\to\infty}\big(P_n^q(S_n(231))\big)^\frac1n$).
Let
\begin{equation}\label{gamma}
\gamma_n=\gamma_n(q)=w_n(q)P_n^q(S_n(312)),
\end{equation}
and define the generating function
$$
G_q(t)=\sum_{n=0}^\infty \gamma_n(q)t^n.
$$
Let $r(q)$ denote the radius of convergence of $G_q$. We note that $r(q)\in[1,\infty)$. Indeed, since
the coefficients of the power series are bounded, we have $r(q)\ge1$, and by
Proposition \ref{nonidentprop} and the
fact that $\lim_{n\to\infty}w_n(q)>0$,
 we have $r(q)<\infty$.
Since $\lim_{n\to\infty}w_n(q)>0$, it follows that
$\limsup_{n\to\infty}\big(P_n^q(S_n(312))\big)^\frac1n=\frac1{r(q)}$, and thus, in light of Proposition \ref{dn1n}, that
\begin{equation}\label{limitingprob}
\lim_{n\to\infty}\big(P_n^q(S_n(312))\big)^\frac1n=
\frac1{r(q)},
 \text{where}\ r(q)\ \text{is the radius of convergence of}\ G_q.
\end{equation}

\begin{proposition}\label{genfuncGprop}
The generating function $G_q$ satisfies the functional equation
\begin{equation}\label{genfuncG}
G_q(t)=\frac1{1-(1-q)tG_q(qt)},\ 0\le t<r(q).
\end{equation}
Furthermore,
\begin{equation}\label{critinfinity}
\lim_{t\to  r(q)}G_q(t)=\infty.
\end{equation}
\end{proposition}
\noindent \bf Remark.\rm\ Note, for example, that the Catalan sequence has generating function
$\frac{1-\sqrt{1-4t}}{2t}$ with radius of convergence $\frac14$, however as $t$ approaches $\frac14$, the generating function remains bounded.

\medskip

Using Proposition \ref{genfuncGprop}, we can prove the following
result.
\begin{theorem}\label{quant}
Let
$$
F(x)=\begin{cases} \frac1{1-x},\thinspace x\in[0,1);\\ \infty,\thinspace x\in\mathbb{R}-[0,1)\end{cases}.
$$
\noindent i. If for some $c\in(0,1]$, one has
\begin{equation}\label{keylower1}
F\big(cF\big(cq\cdots F\big(cq^{N-1}(1+cq^N)\big)\cdots\big)=\infty,\ \text{for some}\ N\in\mathbb{N},
\end{equation}
then
\begin{equation}\label{lowerprob}
\lim_{n\to\infty}\big(P_n^q(S_n(312))\big)^\frac1n=\lim_{n\to\infty}\big(P_n^q(S_n(231))\big)^\frac1n> \frac{1-q}c.
\end{equation}
\noindent ii. If for some $c\in(0,1)$, one has
\begin{equation}\label{keyupper1}
F\big(cF\big(cq\cdots F\big(cq^{N-1}F\big(c\frac{q^N}{1-q}\big)\cdots)<\infty,\ \text{for some}\ N\in \mathbb{N},
\end{equation}
then
\begin{equation}\label{upperprob}
\lim_{n\to\infty}\big(P_n^q(S_n(312))\big)^\frac1n=\lim_{n\to\infty}\big(P_n^q(S_n(231))\big)^\frac1n< \frac{1-q}c.
\end{equation}
In fact, if for some fixed $N$, \eqref{keyupper1} holds with $c$ replaced by $b$, with   $b<c$  arbitrarily close to $c$, then \eqref{upperprob} continues to hold with $c$.
\end{theorem}
\bf\noindent Remark.\rm\
We have restricted $c$ to $(0,1]$ in part (i) because we already know from Proposition \ref{nonidentprop} that
the left hand side of \eqref{lowerprob} is greater or equal to $1-q$. Similarly, we have restricted $c$ to $(0,1)$ in
part (ii) because by Proposition \ref{nonidentprop} and \eqref{upperprob} those are the only possible values of $c$ for which the result could apply.

\begin{figure}\label{fig}
\includegraphics[scale=.75]{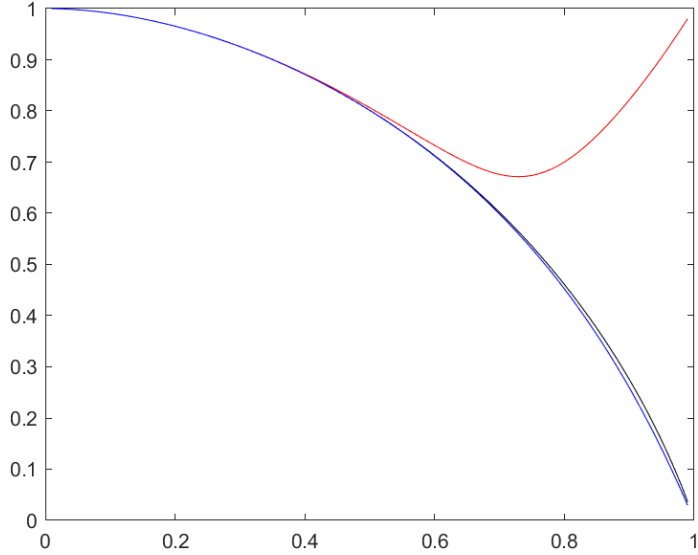}
\caption{Upper bound, lower bound and true value with error $\pm .01$  for $\lim_{n\to\infty}(P_n^q(S_n(312))^\frac1n$ as a function of q }
\end{figure}

The larger one chooses $N$ in \eqref{keylower1} and \eqref{keyupper1}, the better an estimate one obtains.
This works nicely for fixed values of $q$, as we demonstrate in the table below.
However, if we are interested in bounds given as explicit functions of $q$, then we are limited in our choices of $N$.
In \eqref{keylower1}, if we choose $N=2$,
we will get a quadratic inequality for $c$.
Alternatively, since $F$ is increasing on $[0,1)$, the condition
$F\big(cF\big(cq\cdots cq^{N-2}F\big(cq^{N-1}\big)\cdots\big)=\infty$ also implies \eqref{lowerprob}.
Using this condition, we can choose $N=4$ and obtain a quadratic inequality for $c$.
This latter case turns out to give a larger lower bound.
In \eqref{keyupper1}, we obtain a quadratic inequality for $c$ if we choose $N=3$.
However, as the proof of Theorem \ref{quant} shows, the term $\frac{c q^N}{1-q}$
is obtained by making an approximation, and this approximation is only a good one when $\frac{c q^N}{1-q}$ is small.
Because of this, it turns out that the  upper bound using $N=3$ is only reasonably accurate for $q\le.6$. Here is our result.
\begin{theorem}\label{quant2}
Let $q\in(0,1)$. Then
\begin{equation}\label{quant22sided}
\begin{aligned}
&LB(q):=\frac{2q^2(1-q)(1-q^3)}{1-q^4-\sqrt{(1-q^4)^2-4q^2(1-q)(1-q^3)}}<\lim_{n\to\infty}\big(P_n^q(S_n(\tau))\big)^\frac1n<\\
& \frac{2q^2(q^2+1)(1-q)}{1-\sqrt{1-4(1-q)q^2(q^2+1)}}:=UB(q),\ \text{for}\ \tau=312\ \text{and}\ \tau=231.
\end{aligned}
\end{equation}
\end{theorem}
\bf\noindent Remark 1.\rm\
The upper and lower bound functions, $\text{UB}(q)$ and $\text{LB}(q)$, virtually coincide for $q\in(0,.4]$ and differ by less than .01 for $q\in(0,.5]$.
The upper and lower bound functions, $\text{UB}(q)$ and $\text{LB}(q)$, as well as the true value of
 $\lim_{n\to\infty}\big(P_n^q(S_n(312))\big)^\frac1n$ (with error no more than $\pm .01$)
are plotted in  Figure \ref{fig}.
See also the table below.
(We note that the ``true'' values in the table have been obtained by choosing sufficiently large $N$ in \eqref{keylower1} and \eqref{keyupper1} and using MATLAB.)

$$
\begin{aligned}
& q\ \ \ \ \  \  \ \ \ \  \ \ \ \ \ \ \ \  \ \ \ \     .1 \ \ \ \ \ \  \ .2\ \ \ \  \ \ \   .3\ \  \ \ \ \ \  .4\ \ \ \ \ \ \  .5\ \ \ \ \ \ \ .6\ \  \ \ \ \ \  .7\ \  \ \ \ \ \   .8\  \ \ \ \ \ \   .9\\
&\text{UB}(q)\  \ \ \ \   \ \ \  \  \ \ \  \ .991\      \ \ \ \  .966 \ \ \ \   .926 \ \ \ \ .872 \ \ \ \ .806\ \ \ \ .733\ \ \ \ .677\ \ \ \ .700\ \ \ \ .825  \\
& \text{LB}(q) \  \ \ \ \   \ \ \  \  \ \ \  \  .991\  \      \ \ \    .966 \ \ \ \   .926 \ \ \ \ .871 \ \ \ \ .801\ \ \ \ .712\ \ \ \ .599\ \ \ \ .452\  \ \ \ .259  \\
&\text{true value}\ (\pm .01)    \  \  \ \ \ \ \ \ \ \  \ \ \ \ \ \ \  \ \ \ \ \ \ \ \ \   \ \ \  \ \  \  \ \ \ \ \ \ \ \  \ \ \ \ \ \ .716\ \ \ \ .605\ \ \ \ .461\ \ \ \ .275\\
\end{aligned}
$$
\medskip

\bf \noindent Remark 2.\rm\ The upper bound approaches 1 instead of 0 when $q\to1$. Hugo Panzo has shown me the following argument to get the following  rough upper bound (meaningful only for $q>\frac34$) that converges to 0 as $q\to1$:
\begin{equation}\label{Hugo}
\lim_{n\to\infty}\big(P_n^q(S_n(\tau))\big)^\frac1n\le 4(1-q),\ \text{for}\ \tau=312\ \text{and}\ \tau=231.
\end{equation}
This is derived as follows. We have
$$
P_n^q(S_n(312))=\frac1{Z_n(q)}\sum_{\sigma\in S_n(312)}q^{\text{inv}(\sigma)}=\frac1{Z_n(q)}\sum_{k=0}^\infty a_{n,k}q^k,
$$
where $a_{n,k}=|\{\sigma\in S_n(312): \text{inv}(\sigma)=k\}|$ is the number of permutations in $S_n(312)$  that have $k$ inversions.
Since $|S_n(312)|=C_n$, we can rewrite this as
$$
P_n^q(S_n(312))=\frac{C_n}{Z_n(q)}\sum_{k=0}^\infty\frac{a_{n,k}}{C_n}q^k=\frac{C_n}{Z_n(q)}Eq^{A_n},
$$
where $A_n$ is a random variable distributed as the number of inversions occurring in a permutation chosen uniformly at random from $S_n(231)$.
Now \eqref{Hugo} follows from this along with the fact that $\lim_{n\to\infty}C_n^\frac1n=4$ and $\lim_{n\to\infty}Z_n(q)^\frac1n=(1-q)^{-1}$.
\medskip


The rest of the paper is organized as follows. Each of the succeeding sections is devoted to the proof of one result. The results are proved in the following order: Proposition \ref{basicprop},
Proposition \ref{genfuncGprop},
Theorem \ref{quant}, Theorem  \ref{quant2},  Theorem \ref{123thm}, Theorem \ref{132or213}.

\section{Proof of Proposition \ref{basicprop}}\label{proofbasicprop}
\noindent \it Proof of part (i).\rm\ Let $d_n=P_n^q(S_n(312))$.
If $\sigma\in S_n(312)$ and $\sigma_k=1$, then necessarily $\{\sigma_1,\cdots,\sigma_{k-1}\}=[k]\setminus\{1\}$
(and then of course also $\{\sigma_{k+1},\cdots,\sigma_n\}=[n]\setminus[k]$).
Let $A_k\subset S_n$ be the event that $\sigma_k=1$ and $\{\sigma_1,\cdots,\sigma_{k-1}\}=[k]\setminus\{1\}$.
From the online construction,
$$
P_n^q(A_k)=P(X_j\ge1, \text{for}\ j\in[k]-\{1\}; X_{k+l}\le l-1,\ \text{for}\ l\in[n-k]).
$$
Using the fact that
$$
P(X_j\ge1)=\frac{q(1-q^{j-1})}{1-q^j},\ \ \
P(X_{k+l}\le l-1)=\frac{1-q^l}{1-q^{k+l}},\
$$
we obtain
$$
P_n^q(A_k)=q^{k-1}\frac{1-q}{1-q^k}\prod_{l=1}^{n-k}\frac{1-q^l}{1-q^{k+l}}=(1-q)q^{k-1}\frac{w_{k-1}w_{n-k}}{w_n},
$$
where $w_k$ is as in \eqref{w}.
Also from the online construction, it follows that
$$
P_n^q(S_n(312)|A_k)=d_{k-1}d_{n-k}.
$$
Thus,
$$
d_n=\sum_{k=1}^nP_n^q(S_n(312)|A_k)P_n^q(A_k)=(1-q)\sum_{k=1}^nq^{k-1}\frac{w_{k-1}w_{n-k}}{w_n}d_{k-1}d_{n-k},
$$
which is \eqref{basic1}. We leave it to the reader to check that the same formula holds when one works
with 231-avoiding permutations.

\medskip

\noindent \it Proof of part (ii).\rm\
Let $d_n=P_n^q(S_n(213))$.
If $\sigma\in S_n(213)$ and $\sigma_k=1$, then necessarily $\{\sigma_1,\cdots,\sigma_{k-1}\}=\{n-k+2,\cdots, n\}$
(and then of course also $\{\sigma_{k+1},\cdots,\sigma_n\}=\{2,\cdots, n-k+1\}$).
Let $B_k\subset S_n$ be the event that $\sigma_k=1$ and $\{\sigma_1,\cdots,\sigma_{k-1}\}=\{n-k+2,\cdots, n\}$.
From the online construction,
$$
P_n^q(B_k)=P(X_j\le j-2, \text{for}\ j\in[n-k+1]\setminus\{1\}; X_l\ge n-k+1,\ \text{for}\ l\in[n]\setminus[n-k+1]).
$$
Using the fact that
$$
P(X_j\le j-2)=\frac{1-q^{j-1}}{1-q^j};\ \ \ P(X_l\ge n-k+1)=\frac{q^{n-k+1}-q^l}{1-q^l},
$$
we obtain
$$
P_n^q(B_k)=\frac{1-q}{1-q^{n-k+1}}\thinspace q^{(n-k+1)(k-1)}\prod_{l=n-k+2}^n\frac{1-q^{l-n+k-1}}{1-q^l}.
$$
Also from the online construction, it follows that
$$
P_n^q(S_n(213)|B_k)=d_{k-1}d_{n-k}.
$$
Continuing as in the proof of part (i), we obtain \eqref{basic2}.
We leave it to the reader to check that the same formula holds when one works with 132-avoiding permutations.
\hfill $\square$

\section{Proof of Proposition \ref{genfuncGprop}}\label{proofgenfuncGprop}
Using \eqref{basic1} and recalling the definition of $\gamma_n$ in \eqref{gamma}, we have
\begin{equation}\label{gammas}
\gamma_n=(1-q)\sum_{k=1}^n q^{k-1}\gamma_{k-1}\gamma_{n-k}.
\end{equation}
Let $a_n=a_n(q)=q^n\gamma_n$, and define
$$
A_q(t)=\sum_{n=0}^\infty a_n(q)t^{n+1}.
$$
Using  \eqref{gammas}, we have
\begin{equation}\label{GA}
\begin{aligned}
&G_q(t)A_q(t)=\sum_{n=1}^\infty\sum_{k=1}^na_{k-1}\gamma_{n-k}t^n=\sum_{n=1}^\infty\sum_{k=1}^n q^{k-1}\gamma_{k-1}\gamma_{n-k}t^n=\\
&\frac1{1-q}\sum_{n=1}^\infty\gamma_nt^n=\frac1{1-q}(G_q(t)-1).
\end{aligned}
\end{equation}
Since $A_q(t)=tG_q(qt)$, we
conclude from \eqref{GA} that
\begin{equation}\label{Gequ}
G_q(t)=\frac1{1-(1-q)tG_q(qt)}.
\end{equation}

We now prove that $\lim_{t\to r(q)}G_q(t)=\infty$.
Assume to the contrary.
Since $G_q$ is increasing, $\lim_{t\to r(q)}G_q(t)$ exists and from \eqref{Gequ},
$(1-q)r(q)G_q(qr(q))<1$.
 Consequently, for $\delta>0$ sufficiently small, the function
$$
H_q(z)=\begin{cases} G_q(z), |z|<r(q)-\delta;\\ \frac1{1-(1-q)zG_q(qz)},\ r(q)-2\delta<|z|<r(q)+2\delta,\end{cases}
$$
extends analytically the analytic function $G_q(z)$ in $|z|<r(q)$ to the disk $|z|<r(q)+2\delta$.
Consequently, the radius of convergence of $G_q(z)$ must be at least $r(q)+2\delta$, which is a contradiction.
\hfill $\square$

\section{Proof of Theorem \ref{quant}}\label{proofquant}
Recall the definition of  $F$ in the statement of the theorem, and note from \eqref{genfuncG} that
$(1-q)tG_q(qt)<1$ if and only if $t\in[0,r(q))$.
Then \eqref{genfuncG} can be written as
$G_q(t)=F\big((1-q)tG_q(qt)\big)$, for $t\in[0,r(q))$.
In light of  \eqref{critinfinity}, if we   define $G_q(t)=\infty$, for $t\ge r(q)$, then it follows that
\begin{equation}\label{extendedfuncequ}
G_q(t)=F\big((1-q)tG_q(qt)\big),\  \text{for all}\ t\ge0.
\end{equation}

Iterating \eqref{extendedfuncequ}, we conclude that for any $N\in\mathbb{N}$,
\begin{equation}\label{Fiterates}
F\big((1-q)tF\big(q(1-q)t\cdots F\big(q^{N-1}(1-q)t G_q(q^Nt)\big)\cdots\big)<\infty,\ \text{if and only if}\ t\in[0,r(q)).
\end{equation}
Now $G_q(0)=1, G_q'(0)=\gamma_1=1-q$ and $G_q$ is strictly convex, since it is represented as a power series with positive coefficients. Therefore, $G_q(q^Nt)>1+(1-q)q^Nt$, for $t\in(0,r(q))$.
Thus, since $F$ is increasing on $[0,1)$ and is equal to $\infty$ elsewhere, it follows from \eqref{Fiterates} that if
$$
F\big((1-q)tF\big(q(1-q)t\cdots F\big(q^{N-1}(1-q)t(1+(1-q)q^Nt)\big)\cdots\big)=\infty,\
$$
then $r(q)<t$.
The strict inequality $r(q)<t$ follows from the strict inequality $G_q(q^Nt)>1+(1-q)q^Nt$.
Substituting $c=(1-q)t$ above, and using \eqref{limitingprob}, we conclude that if \eqref{keylower1} holds, then \eqref{lowerprob} holds, proving part (i).

We now prove part (ii). Since $\gamma_0=1$ and $\gamma_n<1$, for all $n\ge1$, we have
$G_q(t)<\frac1{1-t}=F(t)$, for $t\in(0,1)$. Thus, since $F$ is  increasing on $[0,1)$ and is equal to $\infty$
elsewhere, it follows from \eqref{Fiterates}  that
\begin{equation}\label{FG}
\text{if}\ \ F\big((1-q)tF\big(q(1-q)t\cdots F\big(q^{N-1}(1-q)t F\big(q^Nt)\big)\cdots\big)<\infty,\ \text{then}\ r(q)>t.
\end{equation}
Substituting $c=(1-q)t$ above, and using \eqref{limitingprob}, we conclude that if  \eqref{keyupper1} holds, then  \eqref{upperprob} holds.

The final statement in part (ii) follows from comparing \eqref{Fiterates} and \eqref{FG} and noting   the strict inequality
$G_q(q^Nt)<F(q^Nt)$, for $t>0$.
Indeed, fix a value of $c\in(0,1)$. If for some fixed value  of $N$,  \eqref{FG} holds for all  $b:=(1-q)t<c$, then by the strict inequality
$G_q(q^Nt)<F(q^Nt)$, for $t\in(0,1)$, it follows that
the left hand side of \eqref{Fiterates} is finite with $(1-q)t=c$, and consequently by \eqref{Fiterates}
and \eqref{limitingprob}, we obtain \eqref{upperprob}.
\hfill $\square$

\section{Proof of Theorem \ref{quant2}}\label{proofquant2}
We begin with the lower bound. As noted in the paragraph preceding Theorem \ref{quant2},
part (i) of Theorem \ref{quant} continues to hold with the condition \eqref{keylower1} replaced by the condition
$F(cF(cq\cdots cq^{N-2}F(cq^{N-1})\cdots)=\infty$.
We use this condition with $N=4$.

Consider the requirement $F\big(cF\big(cqF\big(cq^2F\big(cq^3\big)\big)\big)\big)=\infty$, where we consider  $c\in(0,1)$.
In order for this to occur,  one needs $cF\big(cqF\big(cq^2F\big(cq^3\big)\big)\big)\ge1$.
Since $F(1-c)=\frac1c$, in order for the above  inequality to occur, one needs
$cqF\big(cq^2F\big(cq^3\big)\big)\ge1-c$. Since the range of $F$ is $[1,\infty)$, this second inequality holds automatically if
$\frac{1-c}{cq}\le1$, or equivalently, if
\begin{equation}\label{firstc}
c\ge\frac1{1+q}.
\end{equation}
Otherwise, since
$F(1-\frac{cq}{1-c})=\frac{1-c}{cq}$,
in order for the second  inequality to occur, one needs
$cq^2F\big(cq^3\big)\ge1-\frac{cq}{1-c}=\frac{1-c-cq}{1-c}$.
This third inequality holds automatically if
$\frac{1-c-cq}{(1-c)cq^2}\le 1$, or equivalently, if
\begin{equation}\label{secondc}
c\ge \frac{1+q+q^2-\sqrt{(1+q+q^2)^2-4q^2}}{2q^2}.
\end{equation}
(We don't need to place an upper bound on $c$ because we are restricting from the start to $c\in(0,1)$, and one can see that the inequality
$\frac{1-c-cq}{(1-c)cq^2}\le 1$
holds for $c=1^-$.)
Otherwise,
since
$F(\frac{1-c-cq-cq^2+c^2q^2}{1-c-cq})=\frac{1-c-cq}{(1-c)cq^2}$, in order for the third inequality to hold, one needs
$cq^3\ge\frac{1-c-cq-cq^2+c^2q^2}{1-c-cq}$, or equivalently,
\begin{equation}\label{tocite}
(q^2+q^3+q^4)c^2-(1+q+q^2+q^3)c+1\le0.
\end{equation}
This gives
\begin{equation}\label{thirdc}
c\ge\frac{(1-q^4)-\sqrt{(1-q^4)^2-4q^2(1-q)(1-q^3)}}{2q^2(1-q^3)}.
\end{equation}
(As before, we don't need to place an upper bound on $c$ because we are restricting  to $c\in(0,1)$, and one can see that \eqref{tocite} holds for $c=1$.)

We conclude that \eqref{lowerprob} holds for $c$ satisfying any one of \eqref{firstc}-\eqref{thirdc}. It turns out that the right hand side of \eqref{thirdc} yields the smallest value for $c$, so we choose
$c$ to be equal to the right hand side of \eqref{thirdc}.
 After doing a little algebra, one finds that $\frac{1-q}c$ is given by the left hand side
of \eqref{quant22sided}.
\medskip

We now turn to the upper bound.
We consider \eqref{keyupper1} with $N=3$; that is, we consider the inequality
$F\big(cF\big(cqF\big(cq^2F\big(c\frac{q^3}{1-q}\big)\big)\big)\big)<\infty$, where we consider $c\in(0,1)$.
In order for this to occur,
one needs $cF\big(cqF\big(cq^2F\big(c\frac{q^3}{1-q}\big)\big)\big)<1$.
Since $F(1-c)=\frac1c$, in order for this second  inequality to occur, one needs
$cqF\big(cq^2F\big(c\frac{q^3}{1-q}\big)\big)<1-c$.
If $\frac{1-c}{cq}\le1$, or equivalently, if \eqref{firstc} occurs, then
this second inequality cannot occur. Otherwise, since
 $F(1-\frac{cq}{1-c})=\frac{1-c}{cq}$,
in order for this second inequality to occur, one needs
$cq^2F\big(c\frac{q^3}{1-q}\big)<1-\frac{cq}{1-c}=\frac{1-c-cq}{1-c}$.
 If $\frac{1-c-cq}{(1-c)cq^2}\le1$, or equivalently, if
 \eqref{secondc} occurs, then this third inequality cannot occur. Otherwise,
 since
$F(\frac{1-c-cq-cq^2+c^2q^2}{1-c-cq})=\frac{1-c-cq}{(1-c)cq^2}$, in order for this third inequality to hold, one needs
$\frac{cq^3}{1-q}<\frac{1-c-cq-cq^2+c^2q^2}{1-c-cq}$, or equivalently,
\begin{equation}\label{tociteagain}
(q^2+q^4)c^2-c+1-q>0.
\end{equation}
Thus, we need
\begin{equation}\label{lastc}
c<\frac{1-\sqrt{1-4q^2(1-q)(1+q^2)}}{2q^2(1+q^2)}.
\end{equation}
(We don't need to place a lower bound because we can see that \eqref{tociteagain} holds for $c=0$.)

We conclude that \eqref{keyupper1} holds if $c$ is smaller than each of the three right hand sides, \eqref{firstc}, \eqref{secondc} and \eqref{lastc}.
One can show that the right hand side of \eqref{lastc} is the smallest of the three.
Thus, we choose $c$  to be equal to the right hand side of \eqref{lastc}.  Thus, $\frac{1-q}c$ is given by the right hand side of \eqref{quant22sided}.
\hfill $\square$

\section{Proof of Theorem \ref{123thm}}\label{proof123thm}
Recall that $\{X_j\}_{j=1}^\infty$ are independent and distributed as in \eqref{Xdist}.
We will show that for any $\{i_j\}_{j=1}^m$ satisfying $i_1<\cdots<i_m$,
\begin{equation}\label{finiteinfinite}
\big(\prod_{j=1}^m(1-q^j)\big)\thinspace \frac{q^{\frac12(m-1)m}}{Z_m(q)}\le P(X_{i_1}<X_{i_2}<\cdots< X_{i_m})\le\frac1{\prod_{j=1}^m(1-q^j)}\frac{q^{\frac12(m-1)m}}{Z_m(q)},
\end{equation}
where $Z_m(q)$ is as in \eqref{Zn}.
Using this, we now prove the theorem.
Let
$A_n=\{j\in[n]\}:X_j=j-1\}$ and $B_n=\{j\in[n]:X_j\neq j-1\}$.
It is not hard to see that in the online construction, in order to obtain a 123-avoiding permutation it is necessary
 and sufficient that
$X_{i_1}<X_{i_2}<\cdots<X_{i_m}$, where $m=|B_n|$, $B_n=\{i_l\}_{l=1}^m$ and $i_1<\cdots<i_m$.
Using this with \eqref{finiteinfinite}, and defining $Z_0(q)=1$, we have
\begin{equation}\label{keylower}
\begin{aligned}
&P_n^q(S_n(123))\ge\sum_{k=1}^nP(|A_n|=k)\frac{\prod_{j=1}^{n-k}(1-q^j)}{Z_{n-k}(q)}q^{\frac12(n-k-1)(n-k)}\ge \\
&P(|A_n|\ge m)\frac{\prod_{j=1}^{n-m}(1-q^j)}{Z_{n-m}(q)}\thinspace q^{\frac12(n-m-1)(n-m)},\ \text{for any}\ m\in[n],
\end{aligned}
\end{equation}
and
\begin{equation}\label{keyupper}
\begin{aligned}
&P_n^q(S_n(123))\le\sum_{k=1}^nP(|A_n|=k)\frac1{\prod_{j=1}^{n-k}(1-q^j)}\frac1{Z_{n-k}(q)}q^{\frac12(n-k-1)(n-k)}\le\\
&\sum_{k=1}^nP(|A_n|\ge k)\frac1{\prod_{j=1}^{n-k}(1-q^j)}\frac1{Z_{n-k}(q)}q^{\frac12(n-k-1)(n-k)}.
\end{aligned}
\end{equation}

Now
$$
P(|A_n|\ge m)\ge P(X_j=j-1: j=1,\cdots, m)=\frac{q^{\frac12(m-1)m}}{Z_m(q)}.
$$
Thus,
from \eqref{keylower} we obtain
\begin{equation}\label{keyloweragain}
P_n^q(S_n(123))\ge\frac{\prod_{j=1}^{n-m}(1-q^j)}{Z_{n-m}(q)Z_m(q)}q^{\frac12\big((m-1)m+(n-m-1)(n-m)\big)}, \ \text{for any}\ m\in[n].
\end{equation}
Choosing $m=\lfloor\frac n2\rfloor$ in \eqref{keyloweragain}, we  conclude that
\begin{equation}\label{liminfthm1}
\liminf_{n\to\infty}\big(P_n^q(S_n(123))\big)^{\frac1{n^2}}\ge q^\frac14.
\end{equation}

Also, for some $c=c(q)>0$,
\begin{equation*}
\begin{aligned}
&P(|A_n|\ge k)\le\sum_{1\le i_1<\cdots<i_k\le n}P(X_{i_j}=i_j-1; j=1,\cdots, k)=\\
&\sum_{1\le i_1<\cdots<i_k\le n}\frac{(1-q)^k}{\prod_{j=1}^k(1-q^{i_j})}q^{\sum_{j=1}^k(i_j-1)}\le
\binom nk\frac{(1-q)^k}{\prod_{l=1}^\infty(1-q^l)}\thinspace q^{\frac12(k-1)k}\le\\
&c2^n\thinspace q^{\frac12(k-1)k}.
\end{aligned}
\end{equation*}
Thus, from \eqref{keyupper}, we have for some $c_1=c_1(q)>0$,
\begin{equation}\label{finalupper}
P_n^q(S_n(123))\le c_12^n\sum_{k=1}^nq^{\frac12\big((k-1)k+(n-k-1)(n-k)\big)}\le c_1n2^nq^{\frac{n^2}4+O(n)}.
\end{equation}
From \eqref{finalupper}, we     conclude that
\begin{equation}\label{limsupthm1}
\limsup_{n\to\infty}\big(P_n^q(S_n(123))\big)^{\frac1{n^2}}\le q^\frac14.
\end{equation}
Now  \eqref{123} follows from \eqref{limsupthm1} and \eqref{liminfthm1}.

It remains to prove \eqref{finiteinfinite}.
We begin with the lower bound. We have
\begin{equation}\label{firstlowerthm1}
\begin{aligned}
&\big(\prod_{j=1}^m(1-q^{i_j})\big)P(X_{i_1}< X_{i_2}<\cdots< X_{i_m})=\\
&\sum_{k_1=0}^{i_1-1}\cdots
\sum_{k_{m-1}=k_{m-2}+1}^{i_{m-1}-1}\sum_{k_m=k_{m-1}+1}^{i_m-1}
(1-q)q^{k_1}\cdots(1-q)q^{k_{m-1}}(1-q)q^{k_m}=\\
&\sum_{k_1=0}^{i_1-1}\cdots
\sum_{k_{m-1}=k_{m-2}+1}^{i_{m-1}-1}(1-q)q^{k_1}\cdots(1-q)q^{k_{m-1}}\big(q^{k_{m-1}+1}-q^{i_m}\big)\ge\\
&q(1-q)\sum_{k_1=0}^{i_1-1}\cdots
\sum_{k_{m-1}=k_{m-2}+1}^{i_{m-1}-1}(1-q)q^{k_1}\cdots(1-q)q^{k_{m-2}}(1-q)q^{2k_{m-1}},
\end{aligned}
\end{equation}
where the inequality follows by writing
$$
q^{k_{m-1}+1}-q^{i_m}=q\big(q^{k_{m-1}}-q^{i_m-1}\big)\ge q(1-q)q^{k_{m-1}}.
$$
Continuing with the summation on the right hand side of \eqref{firstlowerthm1}, we have
\begin{equation}\label{secondlowerthm1}
\begin{aligned}
&\sum_{k_1=0}^{i_1-1}\cdots
\sum_{k_{m-1}=k_{m-2}+1}^{i_{m-1}-1}(1-q)q^{k_1}\cdots(1-q)q^{k_{m-2}}(1-q)q^{2k_{m-1}}=\\
&\frac{1-q}{1-q^2}\sum_{k_1=0}^{i_1-1}\cdots
\sum_{k_{m-2}=k_{m-3}+1}^{i_{m-2}-1}(1-q)q^{k_1}\cdots(1-q)q^{k_{m-2}}\big(q^{2k_{m-2}+2}-q^{2i_{m-1}}\big)\ge\\
&(1-q)q^2\sum_{k_1=0}^{i_1-1}\cdots
\sum_{k_{m-2}=k_{m-3}+1}^{i_{m-2}-1}(1-q)q^{k_1}\cdots(1-q)q^{k_{m-3}}(1-q)q^{3k_{m-2}},
\end{aligned}
\end{equation}
where the inequality follows by writing
$$
q^{2k_{m-2}+2}-q^{2i_{m-1}}=q^2(q^{2k_{m-2}}-q^{2i_{m-1}-2})\ge q^2(1-q^2)q^{2k_{m-2}}.
$$
From \eqref{firstlowerthm1} and \eqref{secondlowerthm1} we have
\begin{equation*}
\begin{aligned}
&\big(\prod_{j=1}^m(1-q^{i_j})\big)P(X_{i_1}< X_{i_2}<\cdots< X_{i_m})\ge\\
&(1-q)^2q^{1+2}\sum_{k_1=0}^{i_1-1}\cdots\sum_{k_{m-2}=k_{m-3}+1}^{i_{m-2}-1}(1-q)q^{k_1}\cdots(1-q)q^{k_{m-3}}(1-q)q^{3k_{m-2}}.
\end{aligned}
\end{equation*}
Continuing in the same vein as the above two steps for a total of $m-1$ steps, we obtain,
\begin{equation}\label{thirdlowerthm1}
\begin{aligned}
&\big(\prod_{j=1}^m(1-q^{i_j})\big)P(X_{i_1}< X_{i_2}<\cdots< X_{i_m})\ge\\
&(1-q)^{m-1}q^{1+2+\cdots (m-1)}\sum_{k_1=0}^{i_1-1}(1-q)q^{mk_1}=\\
&(1-q)^{m-1}q^{\frac12(m-1)m}(1-q)\frac{1-q^{mi_1}}{1-q^m}\ge
(1-q)^mq^{\frac12(m-1)m}.
\end{aligned}
\end{equation}
From \eqref{thirdlowerthm1} and \eqref{Zn}, we have
\begin{equation*}
P(X_{i_1}< X_{i_2}<\cdots< X_{i_m})\ge \frac{\prod_{j=1}^m(1-q^j)}{\prod_{j=1}^m(1-q^{i_j})}\thinspace \frac{q^{\frac12(m-1)m}}{Z_m(q)},
\end{equation*}
and the right hand side above is greater than the left hand side of \eqref{finiteinfinite}.

We now turn to the upper bound in \eqref{finiteinfinite}. We have
\begin{equation}\label{firstupperthm1}
\begin{aligned}
&\big(\prod_{j=1}^m(1-q^{i_j})\big)P(X_{i_1}< X_{i_2}<\cdots< X_{i_m})=\\
&\sum_{k_1=0}^{i_1-1}\cdots
\sum_{k_{m-1}=k_{m-2}+1}^{i_{m-1}-1}\sum_{k_m=k_{m-1}+1}^{i_m-1}
(1-q)q^{k_1}\cdots(1-q)q^{k_{m-1}}(1-q)q^{k_m}=\\
&\sum_{k_1=0}^{i_1-1}\cdots
\sum_{k_{m-1}=k_{m-2}+1}^{i_{m-1}-1}(1-q)q^{k_1}\cdots(1-q)q^{k_{m-1}}\big(q^{k_{m-1}+1}-q^{i_m}\big)\le\\
&q\times\sum_{k_1=0}^{i_1-1}\cdots\sum_{k_{m-1}=k_{m-2}+1}^{i_{m-1}-1}(1-q)q^{k_1}\cdots(1-q)q^{k_{m-2}}(1-q)q^{2k_{m-1}}=\\
&q\frac{1-q}{1-q^2}\sum_{k_1=0}^{i_1-1}\cdots
\sum_{k_{m-2}=k_{m-3}+1}^{i_{m-2}-1}(1-q)q^{k_1}\cdots(1-q)q^{k_{m-2}}\big(q^{2k_{m-2}+2}-q^{2i_{m-1}}\big)\le\\
&q\frac{1-q}{1-q^2}q^2\sum_{k_1=0}^{i_1-1}\cdots
\sum_{k_{m-2}=k_{m-3}+1}^{i_{m-2}-1}(1-q)q^{k_1}\cdots(1-q)q^{k_{m-3}}(1-q)q^{3k_{m-2}}.
\end{aligned}
\end{equation}
Continuing in the same vein as the above two steps in \eqref{firstupperthm1} for a total of $m-1$ steps, we obtain
\begin{equation}\label{secondupperthm1}
\begin{aligned}
&\big(\prod_{j=1}^m(1-q^{i_j})\big)P(X_{i_1}< X_{i_2}<\cdots< X_{i_m})\le\\
&q\frac{1-q}{1-q^2}q^2\frac{1-q}{1-q^3}q^3\cdots\frac{1-q}{1-q^{m-1}}q^{m-1}\sum_{k_1=0}^{i_1-1}(1-q)q^{mk_1}\le\\
&q^{\frac12(m-1)m}\frac{(1-q)^m}{\prod_{j=1}^m(1-q^j)}=\frac{q^{\frac12(m-1)m}}{Z_m(q)},
\end{aligned}
\end{equation}
from which the upper bound in \eqref{finiteinfinite} follows.

\hfill $\square$

\section{Proof of Theorem \ref{132or213}}\label{proof132or213}
By Proposition \ref{equals}, it suffices to consider $S_n(213)$.
We use the online construction of a permutation  $\sigma$, distributed as $P_n^q$, using $\{X_j\}_{j=1}^n$, which are independent
 and distributed as in \eqref{Xdist}.

Fix $\alpha\in(\frac12,1)$.
Let $B_{n,\alpha}=\big\{X_j\ge\lfloor n^\alpha\rfloor, \text{for all}\ j\in\{n-\lfloor n^\alpha\rfloor+1,\cdots, n\}\big\}$.
Note that if $B_{n,\alpha}$ does not occur, then  necessarily there exists some
$j\in\{n-\lfloor n^\alpha\rfloor+1, \cdots, n\}$ for which
$|\{i\in\{1,\cdots, n-\lfloor n^\alpha\rfloor\}: \sigma^{-1}_i<\sigma^{-1}_j\}|\ge n-\lfloor 2n^\alpha\rfloor$.
That is, there exists such a $j$ with the property that its position  in the permutation $\sigma$
is to the right of at least $n-\lfloor 2n^\alpha\rfloor$ numbers from among the  numbers
$\{1,\cdots, n-\lfloor n^\alpha\rfloor\}$.  In such case, let
 $$
\{i_1,\cdots, i_m\}= \{i\in\{1,\cdots, n-\lfloor n^\alpha\rfloor\}: \sigma^{-1}_i<\sigma^{-1}_j\},
 $$
  where
 $n-\lfloor 2n^\alpha\rfloor\le m\le n-\lfloor n^\alpha\rfloor$ and 
$1\le i_1<\cdots<i_m\le n-\lfloor n^\alpha\rfloor$.
Let $\sigma'\in S_{n-\lfloor n^\alpha\rfloor}$ be the permutation  obtained by stopping the construction of $\sigma$ after the first
$n-\lfloor n^\alpha\rfloor$ steps.
Under the above scenario, the first $m$ numbers in $\sigma'$ are $\{i_1,\cdots, i_m\}$.
Furthermore, in order for $\sigma$ to belong to $S_n(213)$, it is necessary that 
 the numbers $\{i_1,\cdots, i_m\}$ appear in increasing order
  in $\sigma'$; that is, $\sigma'_1=i_1,\cdots, \sigma'_m=i_m$.

The distribution of $\sigma'$ is the Mallows distribution $P_{n-\lfloor n^\alpha\rfloor}^q$.
From \eqref{mallowsdef}, it follows readily that the  probability  
$$
P_{n-\lfloor n^\alpha\rfloor}^q(\sigma'_1=i_1,\cdots,\sigma'_m=i_m),
$$
as a function of $\{i_1,\cdots, i_m\}$, where
$1\le i_1<\cdots<i_m\le n-\lfloor n^\alpha\rfloor$,
 is maximized when
$\{i_1,\cdots, i_m\}=\{1,2,\cdots, m\}$.
We have 
\begin{equation}\label{altpossibility}
\begin{aligned}
&P_{n-\lfloor n^\alpha\rfloor}^q(\sigma'_1=1,\cdots, \sigma'_m=m)=P(X_1=\cdots=X_m=0)=\\
&\frac{(1-q)^m}{\prod_{j=1}^m(1-q^j)}\le C_1(1-q)^m,
\end{aligned}
\end{equation}
for some constant $C_1$ independent of $n$ and $\alpha$.

From the previous two paragraphs, we deduce that
\begin{equation}\label{almostfinalthm2}
\begin{aligned}
&P_n^q(S_n(213))=P_n^q(B_{n,\alpha}\cap S_n(213))+P_n^q(B_{n,\epsilon}^c\cap S_n(213))\le\\
&P_n^q(B_{n,\alpha})+P_n^q(S_n(213)|B_{n,\epsilon}^c)\le P_n^q(B_{n,\alpha})+C_1(1-q)^{n-\lfloor 2n^\alpha\rfloor}.
\end{aligned}
\end{equation}
From \eqref{Xdist},
\begin{equation}\label{Bnprob}
P_n^q(B_{n,\alpha})=\prod_{j=n-\lfloor n^\alpha\rfloor+1}^n\frac{q^{\lfloor n^\alpha\rfloor}-q^j}{1-q^j}\le
\frac1{\prod_{j=1}^\infty(1-q^j)}\prod_{j=n-\lfloor n^\alpha\rfloor+1}^nq^{\lfloor n^\alpha\rfloor}\le
C_2q^{(\lfloor n^\alpha\rfloor)^2},
\end{equation}
for some constant $C_2$ independent of $n$ and $\alpha$.
From \eqref{almostfinalthm2} and \eqref{Bnprob}, we obtain the upper bound 
$\limsup_{n\to\infty}\big(P_n^q(S_n(213))\big)^\frac1n\le (1-q)$,
 while the lower bound
$\liminf_{n\to\infty}\big(P_n^q(S_n(213))\big)^\frac1n\ge (1-q)$
follows from \eqref{nonident}.\hfill $\square$

\medskip

\noindent \bf Acknowledgment\rm. The author thanks Hugo Panzo for the argument in Remark 2 after Theorem
\ref{quant2}.

\end{document}